\documentclass[12pt,letterpaper]{amsart}

\usepackage{showkeys}
\usepackage{hyperref}
\newcommand{\F}{\ensuremath{\textup{F}}}
\newcommand{\GO}{\ensuremath{{\rm O}}}
\newcommand{\SO}{\ensuremath{{\rm SO}}}

\newcommand{\M}{\ensuremath{{\textup M}}}

\newcommand{\B}{\ensuremath{{\textup B}}}

\newcommand{\HH}{\ensuremath{{\textup H}}}

\newcommand{\T}{\ensuremath{{\textup T}}}
\newcommand{\tr}{\ensuremath{\textup{trace}}}
\newcommand{\V}{\ensuremath{{\textup V}}}
\newcommand{\U}{\ensuremath{{\textup U}}}

\newcommand{\Z}{\ensuremath{\mathbb{Z}}}
\newcommand{\A}{\ensuremath{\mathbb{A}}}
\newcommand{\R}{\ensuremath{\mathbb{R}}}

\newcommand{\Q}{\ensuremath{\mathbb{Q}}}
\newcommand{\Sc}{\ensuremath{\mathcal{S}}}
\newcommand{\OO}{\ensuremath{\mathcal{O}}}

\newcommand{\SLT}{\ensuremath{\widetilde{\textup{SL}}_{2}}}

\newcommand{\rd}[1]{\ensuremath{\textup{ord}( #1 )}}

\numberwithin{equation}{section}
\newtheorem{thm}{Theorem}[section]
\newtheorem{lem}[thm]{Lemma}

\newtheorem{rem}[thm]{Remark}
\newtheorem{definition}[thm]{Definition}

\begin{document}
\title{Average Representation Numbers For Spinor Genera}
\author{Kobi Snitz}
\begin{abstract}
In this paper we establish a formula for the average
of representation numbers of ternary quadratic forms in a 
spinor genus over a totally real number field. 

The significant
fact about the formula is the fact that it is
given in terms of local quantities. Such a formula has
already been established by Kneser \cite{MK} Hsia \cite{JH} and 
Schulze-Pillot
\cite{SP4} by a different method.
\end{abstract}
\maketitle

\section{Introduction}
The problem which motivates this work is that of calculating 
the representation numbers for ternary forms. For a ternary form $Q$
with integral coefficients over a totally real field $\F$ and an integer
$m$ the {\it representation number} $r_{Q}(m)$ is the number of integral
solutions to  $Q(x,y,x)=m$.

Recall the well 
known Seigel formula which gives a weighted sum of
representation numbers $r_{Q_{i}}$ for forms $Q_{1}\ldots Q_{t}$
in a {\it genus} of quadratic forms. Roughly speaking, a genus 
of quadratic forms is a finite collection of forms which are
equivalent over all completions $F_{v}$. The weights in the sum
are given by the size of stabilizers $O(Q_{i})$ of the forms in a
certain group which permutes the forms in a genus. The Seigel formula is
\begin{equation}\label{sg}
\sum_{i=1}^{t}\frac{r_{Q_{i}}(m)}{|\GO(Q_{i})|}=
\left( \sum_{i=1}^{t}\frac{1}{|\GO(Q_{i})|}\right)\prod_{v}\beta_{v}(m)
\end{equation}
where the factors $\beta_{v}$ are {\it local densities} which are
calculated over the completions $\F_{v}$.

Let $\chi=\chi_{\kappa}$ for $\kappa \in \F^{\times}$  be a quadratic
character of $\A^{\times}/\F$ given by the product of {\it Hilbert symbols}
$\chi(a)=\prod_{v}(-\kappa,a_{v})_{v}$. Our formula differs from the
Seigel formula by the presence of a twist by a certain such $\chi$.
\begin{equation} \label{mr}
\sum_{i=1}^{t}\frac{\chi(\nu(h_{i}))}{|O(Q_{i})|}r_{Q_{i}}(m)=
\begin{cases}
C|x|_{\infty}\prod_{v}f_{0,v}(x) 
& \textup{if } m =\kappa x^{2} \\
0
& \textup{if } m \notin \kappa (\F_{v}^{\times})^{2}
\end{cases}
\end{equation}
Where as in the case of the Seigel formula, the right hand side is 
given in terms of locally calculated factors $f_{0,v}$. 
Furthermore,
only a finite collection of the factors $f_{0,v}$ is not identically
$1$ so the product is in fact finite. 
The $h_{i}$ are elements of the group $\HH:=\textup{GSpin}(\V)$ which
permutes the forms $Q_{i}$ and 
$\nu: \HH(\A)\rightarrow \A^{\times}/\A^{\times,2}$ 
is the reduced norm which in this case is also the {\it spinor norm}\cite{JC}.

The composition of $\chi$ with the spinor norm $\nu$, lifts $\chi$ to
a character of $\HH(\A)$. Those $\chi$ for which \eqref{mr} is valid
must be constant on the stabilizers $O(Q_{i})$ of the forms $Q_{i}$. 
The subsets of $h_{1} \ldots h_{t}$ for 
which $\nu(h_{i})=\nu(h_{j})$ form {\it spinor genera} and evidently 
$\chi\circ\nu$ is really a character of the group which permutes the
spinor genera. Thus the characters $\chi$ for which  \eqref{mr} is
valid form a finite group which indexes the spinor genera in a genus \cite{JC}.

By summing over the formulas for all such characters 
we can use finite group Fourier inversion to obtain
the partial sums over the spinor genera.
That is, we can obtain 
 expressions of
the form
\begin{equation}\label{sg-spinor}
\textup{R}(\textup{Spn}(Q),m):=
\sum_{i}\frac{r_{Q_{i}}(m)}{|\GO(Q_{i})|}=
\prod_{v}(\textup{local data})
\end{equation}
 where the sum is over the {\it spinor} genus of $Q$.
\subsection{Acknowledgments}
I would like to thank Rainer Schulze-Pillot
for much helpful discussion about this paper and the 
subject in general. I would also like to thank Stephen Kudla and
Wee Teck Gan for suggesting significant improvements to this paper.

\section{Notation and Basic Definitions}
We will generally follow the notation in \cite{SP1}.
Let $\F$ be a number field and $\OO$ its ring of integers. Let
$\A=\A_{\F}$ be the ring of adeles of $\F$ and denote a completion of
$\F$ by $\F_{v}$. Let $\sigma_{i}:\F\rightarrow\R\;\; i=1 \ldots
\sigma_{d}$ 
be the real embeddings of $\F$. Also, fix a
totally positive quadratic form 
$Q\in\F[x_{1},x_{2},x_{3}]$ 
with integral coefficients.

Let $\V$ be the quadratic space $\F^{3}$ with a
 basis $e_{1},e_{2},e_{3}$ and
let $Q$ also denote the quadratic form on $\V$ defined by 
$Q(a_{1}e_{1}+e_{2}v_{2}+e_{3}v_{3})=Q(a_{1},a_{2},a_{3})$. A full
rank Lattice  $L \subset \V$ with basis $\{l_{1},l_{2},l_{3}\}$
defines a quadratic form $Q_{L}\in\F[x_{1},x_{2},x_{3}]$ by
$Q_{L}(a_{1},a_{2},a_{3})=Q(a_{1}l_{1}+a_{2}l_{2}+a_{3}l_{3})$.
In fact, a lattice $L$ corresponds to a set of integrally equivalent
forms which correspond to different choices of a basis for $L$. 
Evidently, the lattice $L_{0}:=\OO\{e_{1},e_{2},e_{3}\}$ gives rise to
integral equivalence class of the form $Q$.
For an integer $m \in \OO$, the number of points $x \in \OO^{3}$
such that $Q(x)=m$ is called the representation number of $m$ with
respect to $Q$ and is denoted by $r_{Q}(m)$. Notice that the
representation numbers of integrally equivalent lattices are equal 
and therefore lattices are the correct objects to consider when
looking for representation numbers.
\begin{definition} \label{genlat}
Let $L\subset \V$ be a lattice of full rank. The genus of $L$ is the
set of lattices $K\subset \V$ such that the local lattices 
$K_{v}:=K\otimes\OO_{v}$ and $L_{v}:=L\otimes \OO_{v}$ 
are isometric with respect to $Q$ for all local places $v$.
\end{definition}
note the the local equivalence {\it does not} mean that $L$ is
isometric to $K$ over $\F$. 
\begin{definition}
Let $Q_{1}$ and $Q_{2}$ be two quadratic forms corresponding to two
lattices $L_{1}$ and $L_{2}$. $Q_{1}$ and $Q_{2}$ belong to the same 
genus of forms (or simply a genus) if the lattices $L_{1}$ and
$L_{2}$ belong to the same genus.
\end{definition}

A convenient way to describe the genus of $L$ is as the orbit of
the Adelic orthogonal group of $(V,Q)$, $\GO(\A)$ under the action 
$L\rightarrow hL$ for $h =(h_{v})\in \GO(\A)$, where 
$hL\otimes\OO_{v}=h_{v}L_{v}$. The stabilizer of $L_{0}$ under
this action is denoted by $\GO_{\A}(L_{0})$ and $\GO(\A)$ can be
decomposed as the coset space
\begin{equation}\label{decomp}
\GO(\A)=\cup_{i=1}^{t}\GO(\F)h_{i}\GO_{\A}(L_{0}).
\end{equation}
The $h_{i}$ in this decomposition are the
representatives of the lattices in the genus of $L_{0}$ which we
denote by
$L_{i}:=h_{i}\cdot L_{0}$ for $i=1\ldots t$. The corresponding forms 
are denoted by $Q_{i}:=h_{i}\cdot Q$.

For every ternary quadratic form $Q$, there is a quaternion algebra
$Quat(a,b):=\B(\F)$ such that, with the right basis, 
$Q$ is a multiple of the reduced norm 
of the algebra on the space
of trace zero quaternions \cite{VI}. 
We let $\B$ be the algebra that corresponds
to the form $Q$ which we fixed and identify $\V$ with the space of 
trace zero quaternions.
The group $\HH(\A):=\B^{\times}(\A)$ acts on 
$\V(\A)\subset \B(\A)$ by conjugation $h\cdot x=hxh^{-1}$. This action
identifies $\A^{\times}\backslash \HH(\A)$ with $\SO(3)(\A)$.

There are two
possible isomorphism classes of the local algebras
$\B_{v}=\B(\F_{v})$. For almost all $v$ we have $\B_{v}\simeq \M_{2}(\F_{v})$,
the algebra of $2\times 2$ matrices in which case 
$\B_{v}$ is called unramified.
At a finite (even) number of places $v$, the algebra $\B_{v}$ is
isomorphic to a division algebra and we say that $\B_{v}$ is ramified.

\begin{definition}
The spinor norm $\theta: \GO(\F_{v}) \rightarrow  
\F_{v}^{\times}/(\F_{v}^{\times})^{2}$
or
 $\theta: \GO(\F) \rightarrow  
\F^{\times}/(\F^{\times})^{2}$
 \cite{OO} will be used to denote
the spinor norm regardless of the underlying field. With the
identification we made between $\A \backslash \HH(\A) $ and
$\SO(3)(\A)$ the spinor norm for $h_{v}\in \SO(3)(\F_{v})$ is
$Q(h_{v})$.
\end{definition}

\section{The Calculation}
\subsection{Initial Formula}
The basic set up for our calculation is the same as in \cite{SP1} 
and \cite{DSP}.
We use a lift of a one dimensional automorphic 
representation (automorphic character) from
the group $\GO(\A)$ to the group $\SLT(\A)$. 
We will then use the correspondence between automorphic forms and
classical modular forms (see \cite{B}) 
compare Fourier coefficients and extract the arithmetic information.

Let $\chi$ be a quadratic character 
of $\A^{\times}$ given by the product of Hilbert symbols
$\chi(a)=\chi_{\kappa}(a)=\Pi_{v}(a_{v},-\kappa)_{v}$
for some $\kappa \in \F^{\times}$ determined up to a square class.
The automorphic character of $\GO(\A)$ is an extension of the composition 
$\chi\circ Q:\A^{\times}\backslash \HH(\A)\simeq \SO(3)(\A)\rightarrow
\{\pm1\}$. 

Recall that
the representation corresponding to the
character $\chi\circ Q$ is a space of automorphic forms on $\SLT(\A)$, 
each constructed from an integral of a Schwartz function $f \in \Sc(\V(\A))$
multiplied by a theta kernel. The
function $f$ which gives rise to the modular form from which we
extract information is given by
$f:=\Pi_{v}f_{v}$. For every finite $v$ the factor $f_{v}$ is
the characteristic function of $L_{0,v}$ and at the
real places we define $f_{\infty}(x)=\exp(-2\pi\tr(Q(x)))$.

\begin{definition}
Let $\omega$ be the Weil representation of $\SLT(\A)$ on the space
$\Sc(\V(\A))$.
The theta kernel defined for a 
Schwartz function $f\in \Sc(\V(\A))$ a group element $ g \in \SLT(\A)$ 
and an $h\in \GO(\A)$ is 
\[
\theta_{\V}(g,h,f):=\sum_{x\in\V(\F)}\omega(g)f(h^{-1}x)
\]
\end{definition}
With that we can define the theta lift of the character $\chi\circ Q$
\begin{definition}\label{lift}
In the 
lift of the automorphic character $\chi\circ Q$ the Schwartz function
$f$ gives rise to the automorphic form of $g\in\SLT(\A)$
\[
\Theta_{f}(\chi\circ Q)(g):=\int_{\GO(\F)\backslash \GO(\A)}
\chi\circ Q(h)\theta_{\V}(g,h,f)dh
\]
\end{definition}
In order to benefit from the decomposition of $\GO(\A)$ \eqref{decomp}
the character
$\chi$ must be compatible with $Q$ in the following way:
\begin{definition}\label{compat}
The form $Q$ and the character $\chi=\chi_{\kappa}$ 
are said to be compatible if 
\begin{enumerate}
\item $Q$  represents an element 
in the square class of $\kappa$ over $\F$
and
\item
 for every $h\in\GO_{A}(L_{0})$ we have $\chi(Q(h))=1$, i.e.
$\theta(\GO_{\A}(L_{0})\subset \textup{ker}(\chi_{\kappa})$
\end{enumerate}
\end{definition}
\begin{rem}
$\chi$ is compatible with $Q$ if and only if it is compatible with all
  the forms in the genus of $Q$. Therefore we would sometimes say
that $Q$ is compatible with the genus rather than with $Q$.
\end{rem}

The modular form from whose Fourier coefficients contains the
information we seek is the following.
\begin{definition}\label{trace-def}
For $Z =(z_{1},\ldots,z_{d})\in \mathcal{H}^{d}$ 
where $\mathcal{H}$ is the upper half
plane and $d=|\F:\Q|$. For $m \in \F$ the trace $\tr(mz)$ is given by
\[
\tr(mz)=\sum_{i=1}^{d}\sigma_{i}(m)z_{i}
\]
where $\sigma_{i}:\F \rightarrow \R$ are the different embeddings of
$\F$ in $\R$.
Using that, we define the function
\begin{equation*}
\theta_{Q}(Z,L_{0})=\sum_{x \in L_{0}}\exp(2 \pi i \tr(Q(x)Z)).
\end{equation*}
\end{definition}
This function is a modular form of weight $(3/2,\ldots,3/2)$ and
it is not too hard to see that its Fourier expansion is 
\begin{equation} \label{ts}
\theta_{Q}(Z,L_{0})=\sum_{m \in \OO}r_{Q}(m)\exp(2 \pi i \tr(m)Z).
\end{equation}

As in the calculation in \cite{SP1}, 
when the character $\chi$ is compatible with $Q$ we get that for 
$g=(g_{\infty},1\ldots )$
\[
\Theta_{f}(\chi\circ Q)(g)=\sum_{i=1}^{t}
\frac{\chi(Q(h_{i}))}
     {|\GO(h_{i}L_{0})|}
\sum_{x\in h_{i}L_{0}}\omega(g_{\infty})\exp(-2\pi\tr(Q(x)))
\]
where $\GO(h_{i}L_{0})$ is the finite group of integral automorphisms
of the lattice $h_{i}L_{0}$. Next
we have as in the bottom
of page 3 in \cite{SP1} that the corresponding modular form is 
\begin{align} 
\Theta_{f}(\chi\circ Q)(Z)&:=
\sum_{i=1}^{t}\frac{\chi(Q(h_{i}))}{|\GO(Q_{i})|}\theta_{Q_{i}}(Z,L_{0})
\nonumber \\
&=\sum_{i=1}^{t}\frac{\chi(Q(h_{i}))}{|\GO(Q_{i})|}
\sum_{m\in \OO}r_{Q_{i}}(m)\exp(2\pi i \tr(m)Z)\label{bp3}
\end{align}

In the notation of \cite{KS},
equation \eqref{bp3} is the modular form corresponding to the 
automorphic form $I_{\V}(\chi,f)$ where $f=f_{L_{0}}$ defined above.
The main result of \cite{KS} is that 
\[
I_{\V}(\chi,f)=I_{\U}(\mu,f_{0})
\] 
Where $I_{\U}(\mu,f_{0})$ is an automorphic form coming from the lift
of
the character $\mu$ from the one dimensional space $\U$.
Therefore 
we would like to describe $I_{\U}(\mu,f_{0})$ and 
the modular form corresponding to it.
From \cite{KS} we have that 
\begin{equation}
I_{\U}(\mu,f_{0})(g)=\frac{1}{2}\sum_{x\in\F}\omega_{0}(g)f_{0}(x)
\label{sof}
\end{equation}
where $\omega_{0}$ is the Weil representation of $\SLT(\A)$ on the
Schwartz space $\Sc(\U(\A))$. Another result of \cite{KS}
is that the function $f_{0}=\Pi_{v}f_{0,v}$ is defined by the factors
$f_{0,v}\in \Sc(\U(\F_{v}))$ given by 
the local equation
\begin{equation}
f_{0,v}(r)=|r|_{v}\chi_{v}(r)\int_{\T_{x_{0}}(\F_{v})\backslash\HH(\F_{v})}
f(rh^{-1}\cdot x_{0})\chi_{v}(Q_{v}(h))dh
\label{varphi0}
\end{equation}
Where the element $x_{0}$ is any vector in
$\V(\F)$ such that $Q(x_{0})$ is in the square class of $\kappa$, the 
quantity that defines $\chi$ and $\T_{x_{0}}$ is the stabilizer of 
$x_{0}$ in $\HH(\F_{v})$. If no such $x_{0}$ exists then $f_{0,v}=0$
for all $v$.

To simplify \eqref{sof} we need to calculate $f_{0,v}$ for all
$v$. For $x_{0}\in \V$, at almost all places $x_{0,v}\in L_{0,v}$ but
for $r \notin \OO_{v}$ we have $rx_{0}\notin L_{0,v}$. In addition, at
almost all places $\chi_{v}$ and $\B_{v}$ are unramified and 
by lemmas 39  in \cite{KS} (the fundamental lemma) 
we get that $f_{0,v}$
is the characteristic function of $\OO_{v}$ at all such places.
In addition, by lemma 36 in \cite{KS},
at places $v$ where $\B_{v}$ is ramified, if $|x|_{v}$ is small enough
then $f_{0,v}(x)=0$.
Furthermore, by \cite{KS}
we have that for $f_{\infty}(x)=\exp(-2\pi Q(x))$ 
when $Q$ is positive definite
we get $f_{0,\infty}(r)=r\exp(-2\pi\kappa r^{2})$.
putting that together with \eqref{sof} and evaluating at 
$g=(g_{\infty},1,1...)$ we get that
\begin{equation*}
I_{\U}(\mu,f_{0})(g)=\frac{C}{2}
\sum_{x\in\OO}
\prod_{v \in S}f_{0,v}(x)|x|_{\infty}
\omega_{0}(g_{\infty})\exp(-2\pi \tr(\kappa x^{2}))
\end{equation*}
Where $S$ is the finite set of non archimedean places which 
include the dyadic places
and where $f_{0}$ is not the characteristic function of $\OO_{v}$. 
Note also that the expression is
up to scalar $C$ which depends on the normalization of the measure on
$\HH(\A)$.

If we consider the modular form we obtain, from the above automorphic
form we get that 
\begin{equation}
\Theta_{f}(\chi\circ Q)(Z)=\frac{C}{2}\sum_{x\in \OO}|x|_{\infty}
\prod_{v \in S}f_{0,v}(x)\exp(2\pi i \tr(\kappa x^{2})Z)
\label{onedim}
\end{equation}

In \cite{KS} it is proved that $f_{0}$ defined above is even, using
 that and a comparison of \eqref{onedim} with \eqref{bp3} we get
our result:
\begin{thm}
Let $Q$ be an integral quadratic form with lattice $L_{0}$ such
that $Q$ is compatible with $\chi=\chi_{\kappa}$ in the sense 
of definition \ref{compat}. Then for $m \in \OO$ 
\begin{equation}
\sum_{i=1}^{t}\frac{\chi(Q(h_{i}))}{|O(Q_{i})|}r_{Q_{i}}(m)=
\begin{cases}
C|x|_{\infty}\prod_{v\in S}f_{0,v}(x) 
& m =\kappa x^{2} \\
0
&  m \notin \kappa (\F_{v}^{\times})^{2}
\end{cases}
\label{res}
\end{equation}
Where $S$ is the finite set of non archimedean places 
which includes the diadic places
and places where the space $\V_{v}$ is
not split or the character $\chi_{\kappa}$ is ramified.
\end{thm}
\subsection{Averages Over Spinor Genera}
It is possible to use formula \eqref{res} 
and others like it to obtain
information about the average number of representations over a spinor
genus. Recall that two forms $Q_{i}$ and $Q_{j}$ are in the same spinor
genus if and only if the rotation relating them is in the kernel of
the spinor norms, i.e.
$Q(h_{i,v}h_{j,v})\in (\F_{v}^{\times})^{2}$ for
all $v$. This means that $\chi$ is constant on each spinor genera and
the sum in \eqref{res} breaks up into sums over the different spinor
genera.
\begin{equation}\label{res-spn}
\sum_{j=1}^{l}\chi(Q(h'_{j}))
\textup{R}(\textup{Spn}(Q'_{j}),m)=
\begin{cases}
C|x|_{\infty}\prod_{v\in S}f_{0,v}(x) 
& m =\kappa x^{2} \\
0
&  m \notin \kappa (\F_{v}^{\times})^{2}
\end{cases}
\end{equation}

where we have relabeled $h'_{j}$ as the representatives of the spinor
genera $Q'_{j}:=h'_{j}\cdot Q\;\;\;j=1 \ldots l$.

The number of spinor genera in a genus corresponds to the number
of elements of the Abelian 
group $G=R/ST$ defined in \cite{JC} p. 209. 
It is immediate from the definition of $G$ that
its characters correspond to characters which
satisfy the second condition of compatibility with $Q$ in definition 
\ref{compat}.
Furthermore, The characters of $G$ are quadratic and therefore they
all correspond to $\chi_{\kappa_{j}}$ for some $\kappa_{j}\in \F^{\times}$. 

From the definition of $f_{0,v}$ we can see that a character $\chi$ 
of  $G$ is compatible with the form $Q$
if and only if the right hand side of equation \eqref{res-spn} is not
identically $0$. Whether or not $\chi$ is compatible with the genus,
there will be $l$ equations of the form of \eqref{res-spn} with
$\chi$ being one of $\chi_{\kappa_{1}}\ldots \chi_{\kappa_{l}}$.For
those $\chi_{k_{j}}$ which are not compatible, the right hand side of 
\eqref{res-spn} will be zero but in any case the
orthogonality relations of the characters of $G$ allow us to solve for
each one of the $\textup{R}(\textup{Spn}(Q'_{j}),m)$ for all
$m\ \in \OO$ and $j=1 \ldots l$. For the sake of clarity is should be
said that the question of compatibility is only relevant because it
tells us in which cases the right hand side of \eqref{res-spn} is
already known.

\begin{lem}\label{when-compat}
When $\F=\Q$ the character $\chi=\chi_{\kappa_{j}}$ as above is 
always compatible with the form $Q$. When $\F\neq \Q$, checking 
the condition of compatibility reduces to checking for compatibility 
at the diadic places.
\end{lem}

\begin{proof}
It is a fact that $Q$ 
represents an element from the square class of $\kappa_{j}$
if and only if $\kappa_{j}$ is locally represented
at all places. The question of local representability of $\kappa_{j}$
almost always follows from the second condition of 
compatibility of $Q$ and
$\chi_{\kappa_{j}}$.
The spinor norms of $\SO(3)(\F_{v})$ are given in 
\cite{EH} \cite{FX} and \cite{MK} and they allow us to verify that at 
least at the odd primes, when 
$\theta(\SO_{\F_{v}}(L_{0,v}))\subset \textup{ker}(\chi_{\kappa})$
then $\kappa$ is represented by $L_{0,v}$. 

When we restrict to $\F=\Q$ then we can also verify the above
condition at $\Q_{2}$ when $a$ or $b$ which define the quaternion
algebra $\textup{Quat}(a,b)$ 
are units. Therefore, at least over $\Q$, when
2 does not divide $a$ or $b$, then
all the characters of $G$ are compatible with the genus.
\end{proof}

It is the author's intention 
to carry the calculation
further in a follow up paper. In particular, it is possible to
transform the orbital integrals defining $f_{0,v}$ into line
integrals and simplify them significantly. It is also possible to 
use the formulas in \cite{TY} to simplify the local densities 
$\beta_{v}$.

In some situations (when the spinor genera consist of a single class)
it is possible to use this result even to obtain
formulas for representation numbers of individual forms.

\section{Example}
We start with the form 
\begin{equation}\label{q1}
Q_{1}(x,y,z)=4x^{2}+16y^{2}+64z^{2}.
\end{equation} 
Using the Kneser neighbor method \cite{KNE} implemented by a computer
program of Schulze-Pillot \cite{SP3} we get that the genus
of this form has one additional form
\begin{equation}\label{q2}
Q_{2}(x,y,z)=20x^{2}+16xy +16y^{2}+16z^{2}.
\end{equation}
and $|O(Q_{1})|=|O(Q_{2})|=8$. In order to find a character compatible
with the form we need to consider the set 
$\theta(\SO_{\Q_{v}}(L_{0,v}))$ for
every $v$. These are 
described in \cite{EH} and \cite{JC}. In
particular, lemma 3.3 p. 208 in \cite{JC} says that for the form $Q$,
at all odd primes $v$, 
$\theta(\SO_{\Q_{v}}(L_{0,v}))=\OO_{v}^{\times}(\Q_{v}^{\times})^{2}$. 
This means that in order for a character $\chi_{\kappa}$ to be
compatible with the form $Q$ it must be trivial on 
$\OO_{v}^{\times}(\Q_{v}^{\times})^{2}$ for all odd $v$. That is,
$(-k,a)_{v}=1$ iff $\rd{a}$ is even, i.e. $\rd{-\kappa}$ is also even.
The spinor norms of the stabilizers of the lattice
over $\Q_{2}$ is given by theorem 2.7 in \cite{EH} and in
our case it is equal to the kernel of $\chi_{1}$. Therefore the
character $\chi=\chi_{4}$ is trivial on the stabilizer of the 
lattice and since $Q_{1}$ represents $4$ the character is compatible
with $Q_{1}$.
 To further evaluate equation \eqref{res} we use the 'fundamental
 lemma' of \cite{KS} which says that whenever $L_{0,v}$ is a split
 lattice and the character $\chi_{\kappa}$ is unramified then
 $f_{0,v}$ is the characteristic function of $\OO_{v}$. The lattice of
 the form $Q_{1}$ is split at all odd primes and the character
 $\chi_{4}$ is unramified at all odd primes.
Therefore for the genus of two forms $Q_{1}$ and 
$Q_{2}$, equation \eqref{res} becomes.
\begin{equation}
r_{Q_{1}}(m)-r_{Q_{2}}(m)=
\begin{cases}
8C\sqrt{m}f_{0,2}(\sqrt{m})
& m \in (\Z^{\times})^{2}\\
0
&  m \notin (\Z^{\times})^{2}
\end{cases}
\label{ressimp}
\end{equation}
where 
\begin{equation}\label{f02}
f_{0,2}(x)=|x|_{2}(-1,x)_{2}\int_{\T_{x_{0}}(\Q_{2})\backslash\HH(\Q_{1})}
\textup{Char}_{L_{0,2}}(xh^{-1}\cdot x_{0})\chi_{2}(Q_{2}(h))dh
\end{equation}
We can simplify the evaluation of equation \eqref{ressimp} by noticing
the dependence of \eqref{f02} on $x$. It is not hard to see that if
$n$ is an odd number then
\begin{equation}\label{trans} 
f_{0,2}(nx)=(-1,n)_{2}f_{0,2}(x)
\end{equation}
Next, it is clear that neither of the forms $Q_{1}$ or $Q_{2}$ can
represent odd numbers.
These observations give us the right hand side of 
\eqref{ressimp} for all integers $m=2^{j}n$ for $n$ odd all odd powers $j$ as
well as $j=0,2$. For the remaining even powers of $j$ we use lemma 5.11
in \cite{KS} which says that at a prime $v$ where $\chi$ is non trivial and 
the lattice $L_{0,v}$ is non split (such is the case at $v=2$ in this
example) then $f_{0,2}(2^{j}n)$ vanishes for large enough $j$. The
value $j$ such that $f_{0,2}(2^{l}n)= 0$ for all $l \geq j$ can be determined 
following the proof of lemma 5.11 
in \cite{KS}. The condition for vanishing is that
if $v\in V(\Q_{2})$ and $Q_{1}(v)=2^{j}n$ then $v \in L_{0,2}$. In our
case this amounts to verifying that the solutions to 
$4 \alpha^{2} + 16\beta^{2} + 64\gamma^{2}=2^{j}n$ must be 2-adic integers 
for all $j\geq 6$.
The remaining cases of $m=4,16$ are calculated by directly calculating
$r_{Q_{2}}(4)=0,\;r_{Q_{2}}(16)=4,\;r_{Q_{1}}(4)=2,\;
;r_{Q_{1}}(16)=4$. Using equation \eqref{ressimp}
this gives that $8Cf_{0,2}(2)=1$ and $8Cf_{0,2}(4)=0$
and we finally get that
\begin{equation}\label{resfin}
r_{Q_{1}}(m)-r_{Q_{2}}(m)=
\begin{cases}
2n(-1,2n)_{2} & m=4n^{2} \textup{ for odd } n \\
0 & \textup{ otherwise }
\end{cases}
\end{equation} 

We can combine \eqref{resfin} and the Seigel formula \eqref{sg}
to find expressions for $r_{Q_{1}}(m)$ and $r_{Q_{2}}(m)$ which depend
only on local calculations.
The Seigel formula for the genus in this example is
\begin{equation}\label{sgex}
r_{Q_{1}}(m) + r_{Q_{2}}(m)=2\prod_{v}\beta_{v}(m)
\end{equation} 
Where the product is taken over all primes.
The $\beta_{v}(m)$ are the local densities defined by
Seigel and they 
were calculated explicitly in \cite{TY}.
Thus we have obtained the following
\begin{equation}\label{r1}
r_{Q_{1}}(m) =
\begin{cases}
n(-1,2n)_{2} + 2\prod_{v}\beta_{v}(m) & 
m=4n^{2} \textup{ for odd } n \\
2\prod_{v}\beta_{v}(m) & \textup{ otherwise }
\end{cases}
\end{equation}
and
\begin{equation}\label{r2}
r_{Q_{2}}(m) =
\begin{cases}
2\prod_{v}\beta_{v}(m) - n(-1,2n)_{2} & m=4n^{2} \textup{ for odd } n \\
2\prod_{v}\beta_{v}(m) & \textup{ otherwise }
\end{cases}
\end{equation}
These formulas give the representation numbers in terms of local
factors but we can do a little better with a recursion formula and
reduce the calculation to a finite number of factors.
First of all, the product of the archimedean
factors is given as $C\sqrt{m}$ where $C$ is some constant which does
not depend on $m$. Next, 
In the notation in \cite{TY} $\beta_{v}(m)=\alpha_{v}(m,S)$ where $S$
is the gram matrix of the form $Q_{1}$.
It follows from the definition
of the local densities that if $u\in \Z^{\times}_{v}$ then
$\beta_{v}(u^{2}m)=\beta_{v}(m)$. Therefor we obtain that if 
$d=p_{1}^{2k_{1}}p_{2}^{2k_{2}}\ldots p_{r}^{2k_{r}}$ for distinct primes 
$p_{1}\ldots p_{r}$ then
\begin{align}
\prod_{v}\beta_{v}(dm) &= 
\sqrt{d}\prod_{i=1}^{r}\frac{\beta_{p_{i}}(nm)}
{\beta_{p_{i}}(m)}\prod_{v}\beta_{v}(m)\\
&= \sqrt{d}\prod_{i=1}^{r}\frac{\beta_{p_{i}}(dm)}{\beta_{p_{i}}(m)}
\left[ \frac{r_{Q_{1}}(m) +r_{Q_{2}}(m)}{2} \right]
\end{align}

and our final formulas are
\begin{equation}\label{r1fin}
r_{Q_{1}}(dm) =
\begin{cases}
n(-1,2n)_{2} +  
\sqrt{d}\prod_{i=1}^{r}\frac{\beta_{p_{i}}(dm)}{\beta_{p_{i}}(m)}
R(m)& 
dm=4n^{2} \textup{ for odd } n \\
\sqrt{d}\prod_{i=1}^{r}\frac{\beta_{p_{i}}(dm)}{\beta_{p_{i}}(m)}
R(m) & \textup{ otherwise }
\end{cases}
\end{equation}
and
\begin{equation}\label{r2fin}
r_{Q_{2}}(dm) =
\begin{cases}
\sqrt{d}\prod_{i=1}^{r}\frac{\beta_{p_{i}}(dm)}{\beta_{p_{i}}(m)}
R(m) - n(-1,2n)_{2} & m=4n^{2} \textup{ for odd } n \\
\sqrt{d}\prod_{i=1}^{r}\frac{\beta_{p_{i}}(dm)}{\beta_{p_{i}}(m)}
R(m)
 & \textup{ otherwise }
\end{cases}
\end{equation}
where $R(m)=r_{Q_{1}}(m) +r_{Q_{2}}(m)$
\bibliography{kbib}
\bibliographystyle{amsplain}

\end{document}